\documentclass{article}

\usepackage{geometry}
\usepackage{graphicx}
\usepackage{multirow}
\usepackage{amsmath,amssymb,amsfonts}
\usepackage{amsthm}
\usepackage{mathrsfs}
\usepackage[title]{appendix}
\usepackage{xcolor}
\usepackage{textcomp}
\usepackage{manyfoot}
\usepackage{booktabs}
\usepackage{algorithm}
\usepackage{algorithmicx}
\usepackage{algpseudocode}
\usepackage{listings}
\usepackage{enumerate}
\usepackage{tikz}

\newcommand{\des}{\operatorname{des}}

\newcommand{\Min}{\operatorname{Min}}
\newcommand{\Max}{\operatorname{Max}}
\newcommand{\Lext}{\mathcal L}
\newcommand{\R}{\mathcal R}

\theoremstyle{plain}
\newtheorem{theorem}{Theorem}
\newtheorem{proposition}[theorem]{Proposition}
\newtheorem{corollary}[theorem]{Corollary}
\newtheorem{conjecture}[theorem]{Conjecture}
\newtheorem{lemma}[theorem]{Lemma}

\theoremstyle{remark}

\theoremstyle{definition}

\title{Balanced Simplicial Complexes for the $\gamma$-Vectors of
	Zig-Zag Eulerian Polynomials}

\author{%
	Wenle Shi and Lijie Wang\thanks{%
		Corresponding author: Lijie Wang.\\
		Email addresses: \texttt{shi-wenle@hotmail.com},
		\texttt{lijiewang26@hotmail.com}.%
	}\\[0.5em]
	\small School of Mathematical Sciences, Dalian University of Technology\\
	\small Dalian 116024, P. R. China
}

\date{}

\begin{document}
	
	\maketitle
	
	\begin{abstract}
		In this paper, we construct a balanced simplicial complex whose face numbers give the $\gamma$-coefficients of the descent polynomial of any canonically labeled rank-one poset.
		As an application, we prove Petersen and Zhuang's conjecture that the $\gamma$-vector of the normalized zig-zag Eulerian polynomial is the $f$-vector of a balanced simplicial complex.
		We also obtain analogous results for equatorial spheres of rank-one posets and for the Ehrhart $h^*$-polynomials of their order and chain polytopes.
		
		\medskip
		\noindent\textbf{Keywords:} Zig-zag Eulerian polynomials, $\gamma$-vectors,
		balanced simplicial complexes, Frankl--F\"uredi--Kalai inequalities
		
		\smallskip
		\noindent\textbf{MSC Classification:} 05E45, 05A05
	\end{abstract}
	
	\section{Introduction}

The zig-zag poset (fence) $P_n$ has elements $p_1,\ldots,p_n$ and covering relations
\begin{align*}
p_1<p_2>p_3<p_4>\cdots.
\end{align*}
A labeling $\omega:P\to\mathbb Z$ of a finite poset $P$ is \emph{natural} if it is injective and $x<_P y$ implies $\omega(x)<\omega(y)$.
For a finite poset $P$ with $n=|P|$, a linear extension is an ordering $x_1,\ldots,x_n$ of the elements of $P$ that respects all order relations.
Let $\Lext(P,\omega)$ be the set of label words $\omega(x_1)\cdots\omega(x_n)$ of the linear extensions of $P$. 
For a label word $\pi=\pi_1\cdots\pi_n$, let $\des(\pi)$ denote the number of indices $i$ with $1\le i<n$ such that $\pi_i>\pi_{i+1}$.

Following Petersen and Zhuang~\cite[Section~2.2]{PZ}, we define the \emph{zig-zag Eulerian polynomial} by
\begin{align*}
Z_n(t)=\sum_{\pi\in\Lext(P_n,\omega)}t^{1+\des(\pi)},
\end{align*}
where $\omega$ is a natural labeling of $P_n$. 
Petersen and Zhuang~\cite[Theorem~2.5]{PZ} proved that, for $n\ge2$,
\begin{align}\label{eq:gamma}
U_n(t):=\frac{Z_n(t)}{t}
=\sum_{k=0}^{\lfloor\frac{n-2}{2}\rfloor}\gamma_{n,k}
t^k(1+t)^{n-2-2k},
\end{align}
where the coefficients $\gamma_{n,k}$ are nonnegative integers.
The vector
\begin{align*}
\gamma_n=(\gamma_{n,0},\ldots,\gamma_{n,\lfloor(n-2)/2\rfloor})
\end{align*}
is called the \emph{$\gamma$-vector} of $U_n(t)$.
For $n=1$, we have $U_1(t)=1$ and set $\gamma_1=(1)$.
The nonnegative expansion in~\eqref{eq:gamma} implies that the coefficients of $U_n(t)$ are symmetric and unimodal.
Alexandersson~\cite[Theorem 7.1]{Alex} proved by a differential-operator argument that $U_n(t)$ has only simple negative roots for $n\geq 3$. 
Torres~\cite[Theorem~1.1]{Torres} also obtained this result using a transfer matrix and complex analysis.
By Newton's inequalities, the coefficients of $U_n(t)$ are log-concave. 
Alexandersson and Jal~\cite{AlexJal} proved that the $P$-Eulerian polynomials of naturally labeled posets of width two are ultra-log-concave.

For a simplicial complex $\Delta$ of dimension $d-1$, let $f_{k-1}(\Delta)$ denote the number of faces with $k$ vertices. 
In particular, $f_{-1}(\Delta)=1$, and $f_{i}(\Delta)=0$ for $i>\dim \Delta$.
The \emph{face vector} ($f$-vector)  of $\Delta$ is
\begin{align*}
f(\Delta)=\bigl(f_{-1}(\Delta),\ldots,f_{d-1}(\Delta)\bigr).
\end{align*}
For a simplicial complex $\Delta$, write $\dim\Delta=\max_{F\in\Delta}(|F|-1)$.
A simplicial complex $\Delta$ with $\dim\Delta=d-1$ is \emph{balanced} if its vertices admit a proper coloring with $d$ colors, meaning that distinct vertices in any face have distinct colors.
The definition applies to both pure and nonpure complexes.
The complex $\{\varnothing\}$ has $\dim\{\varnothing\}=-1$ and is balanced with zero colors.

Nevo and Petersen~\cite{NP} constructed explicit flag complexes whose $f$-vectors are the $\gamma$-vectors of certain flag homology spheres. 
They also showed that these $\gamma$-vectors satisfy the Kruskal--Katona and Frankl--F\"uredi--Kalai inequalities.
Nevo, Petersen, and Tenner~\cite{NPT} proved that the $\gamma$-vector of the barycentric subdivision of a simplicial sphere is the $f$-vector of a balanced simplicial complex.

A nonempty finite poset $P$ has \emph{rank one} if every maximal chain of $P$ contains exactly two elements. 
A \emph{canonical labeling} of a rank-one poset $P$ is an injective map $\omega:P\to\mathbb Z$ such that $\omega(x)<0$ for every minimal element $x$ and $\omega(y)>0$ for every maximal element $y$.
If $x<_P y$, then $x$ is minimal and $y$ is maximal, so $\omega(x)<0<\omega(y)$. Thus every canonical labeling is natural.

The following theorem gives our construction for every canonically labeled rank-one poset.

\begin{theorem}\label{thm:general}
For every canonically labeled rank-one poset $(P,\omega)$ on $n\ge2$ elements, there is an explicit simplicial complex $\Delta(P,\omega)$ with a proper $\lfloor(n-2)/2\rfloor$-coloring such that
\begin{align*}
\sum_{\pi\in\Lext(P,\omega)}t^{\des(\pi)}
=
\sum_{k=0}^{\lfloor(n-2)/2\rfloor}
f_{k-1}(\Delta(P,\omega))\,t^k(1+t)^{n-2-2k}.
\end{align*}
If both rank levels have even size and every minimal element is below every maximal element, then $\dim\Delta(P,\omega)=\lfloor(n-2)/2\rfloor-2$.
Otherwise, $\dim\Delta(P,\omega)=\lfloor(n-2)/2\rfloor-1$.
\end{theorem}

Frankl, F\"uredi, and Kalai~\cite[Theorem~1.2]{FFK} characterized the $f$-vectors of $d$-colorable simplicial complexes by upper bounds on $f_{i+1}$ in terms of $f_i$, known as the Frankl--F\"uredi--Kalai inequalities.
Petersen and Zhuang~\cite[Conjecture~6.1(c)]{PZ} proposed the following conjecture.

\begin{conjecture}\label{conj:fence}
For every $n\ge1$, the $\gamma$-vector $\gamma_n$ of $U_n(t)$ satisfies the Frankl--F\"uredi--Kalai inequalities. 
Equivalently, $\gamma_n$ is the $f$-vector of a balanced simplicial complex.
\end{conjecture}

The following result proves Conjecture~\ref{conj:fence}, giving an explicit construction and a dimension formula.

\begin{theorem}\label{cor:fence}
For every $n\ge1$, there is an explicitly defined balanced simplicial complex $\Delta_n$ whose $f$-vector is $\gamma_n$. Moreover,
\begin{align*}
\dim\Delta_n
=\max\bigl\{-1,\lfloor(n-2)/2\rfloor-1\bigr\}.
\end{align*}
Consequently, $\gamma_n$ satisfies the Frankl--F\"uredi--Kalai inequalities.
\end{theorem}

The paper is organized as follows.
In Section~2, we give a combinatorial interpretation of the $\gamma$-coefficients and construct a simplicial complex whose face numbers equal these coefficients.
In Section~3, we prove Theorems~\ref{thm:general} and~\ref{cor:fence}.
Section~\ref{sec:applications} gives further applications.

\section{The ideal-chain complex and face numbers}

Throughout the paper, all posets are finite.
In a rank-one poset, the minimal and maximal elements form two disjoint rank levels, and every strict relation is of the form a minimal element to a maximal element. 

Let $\pi=\pi_1\cdots\pi_n\in\Lext(P,\omega)$.
An index $i$ with $1\le i<n$ is a \emph{descent position} if $\pi_i>\pi_{i+1}$. Set $\pi_0=\pi_{n+1}=0$.
An index $i$ with $1\le i\le n$ is a \emph{double descent} if $\pi_{i-1}>\pi_i>\pi_{i+1}$.

For completeness, we recall why the choice of natural labeling does not change the descent polynomial. 
Let $n=| P|$ and let $\omega$ be a natural labeling of $P$. 
If $\Omega_P(q)$ denotes the number of order-preserving maps $P\to\{1,\ldots,q\}$, the classical $P$-partition identity~\cite[Eq.~(1)]{PZ} is
\begin{equation}\label{eq:Ppartition}
\sum_{q\ge1}\Omega_P(q)t^q
=\frac{t\sum_{\pi\in\Lext(P,\omega)}t^{\des(\pi)}}{(1-t)^{n+1}}.
\end{equation}
Its left-hand side depends only on $P$, so the descent polynomial on the right is independent of the natural labeling. 

A cover relation has sign $+1$ if the labels increase along the relation and sign $-1$ if they decrease.
A labeled poset is \emph{sign-graded} if the sum of these signs is the same along every maximal chain.
The common sum is called the \emph{sign-rank}.

For sign-graded posets, we use canonical labelings in the sense of Br\"and\'en~\cite[Section~6]{Branden}.
For $\pi\in\mathcal L(P,\omega)$, write $\operatorname{Orb}(\pi)$ for the orbit of $\pi$ under the $\mathbb Z_2^P$-action generated by the involutions in~\cite[Definition~6.2]{Branden}.

\begin{lemma}[{\cite[Theorem~6.3]{Branden}}]
\label{thm:branden}
Let $(P,\omega)$ be a canonically labeled sign-graded poset of sign-rank $r$, and let $n=|P|$.
For every $\pi\in\mathcal L(P,\omega)$, the orbit $\operatorname{Orb}(\pi)$ contains a unique representative $\widehat{\pi}$ such that $0\widehat{\pi}0$ has no double descents. 
Moreover,
\begin{align*}
\sum_{\sigma\in\operatorname{Orb}(\pi)}
t^{\des(\sigma)}
=t^{\des(\widehat{\pi})}
(1+t)^{n-r-1-2\des(\widehat{\pi})}.
\end{align*}
\end{lemma}

Let $(P,\omega)$ be a canonically labeled rank-one poset on $n$ elements. 
Every maximal chain consists of a single cover relation of sign $+1$, so $(P,\omega)$ is sign-graded with sign-rank $1$. 
The generalized rank function is $0$ on minimal elements and $1$ on maximal elements. 
Since $\omega$ assigns negative labels to elements of rank $0$ and positive labels to elements of rank $1$, it satisfies Br\"and\'en's canonical labeling conditions. 
Applying Lemma~\ref{thm:branden} with $r=1$, we obtain
\begin{equation}\label{eq:rank-one-orbit}
\sum_{\sigma\in\operatorname{Orb}(\pi)}
t^{\des(\sigma)}
=t^{\des(\widehat{\pi})}
(1+t)^{n-2-2\des(\widehat{\pi})}.
\end{equation}

For each integer $k\ge0$, define
\begin{equation}\label{eq:Rk}
\mathcal R_k(P,\omega)
=\left\{
\pi\in\mathcal L(P,\omega):
\des(\pi)=k,\ 
0\pi0\text{ has no double descents}
\right\}.
\end{equation}

\begin{proposition}\label{prop:gamma-count}
Let $(P,\omega)$ be a canonically labeled rank-one poset on $n\ge2$ elements.
Then
\begin{equation}\label{eq:orbit}
\sum_{\pi\in\Lext(P,\omega)}t^{\des(\pi)}
=\sum_{k=0}^{\lfloor(n-2)/2\rfloor}|\R_k(P,\omega)|
t^k(1+t)^{n-2-2k}.
\end{equation}
\end{proposition}

\begin{proof}
By Lemma~\ref{thm:branden}, each orbit of the $\mathbb Z_2^P$-action on $\mathcal L(P,\omega)$ contains a unique representative $\widehat{\pi}$ such that $0\widehat{\pi}0$ has no double descents.
By the definition in~\eqref{eq:Rk}, this representative belongs to $\mathcal R_k(P,\omega)$ precisely when $\des(\widehat{\pi})=k$.

For $\widehat{\pi}\in\mathcal R_k(P,\omega)$, the orbit formula~\eqref{eq:rank-one-orbit} states that
\begin{align*}
\sum_{\sigma\in\operatorname{Orb}(\widehat{\pi})}
t^{\des(\sigma)}
=t^k(1+t)^{n-2-2k}.
\end{align*}
Since the left-hand side is a polynomial and $t^k$ is not divisible by $1+t$, we must have $n-2-2k\ge0$.
Hence every orbit representative belongs to $\mathcal R_k(P,\omega)$ for some $0\le k\le \lfloor(n-2)/2\rfloor$.

Distinct representatives belong to distinct orbits, and every linear extension belongs to exactly one orbit.
Therefore, summing over the orbits indexed by these representatives counts every element of $\mathcal L(P,\omega)$ exactly once.
Using~\eqref{eq:rank-one-orbit}, we obtain
\begin{align*}
\sum_{\pi\in\mathcal L(P,\omega)}t^{\des(\pi)}
=
\sum_{k=0}^{\lfloor(n-2)/2\rfloor}
\sum_{\widehat{\pi}\in\mathcal R_k(P,\omega)}
\sum_{\sigma\in\operatorname{Orb}(\widehat{\pi})}
t^{\des(\sigma)}
=
\sum_{k=0}^{\lfloor(n-2)/2\rfloor}
|\mathcal R_k(P,\omega)|
t^k(1+t)^{n-2-2k}.
\end{align*}
This completes the proof.
\end{proof}

The polynomials $t^k(1+t)^{n-2-2k}$ for $0\le k\le\lfloor(n-2)/2\rfloor$ are linearly independent because their lowest nonzero degrees are distinct. 
Let $\gamma_k(P)$ denote the coefficient of $t^k(1+t)^{n-2-2k}$ in \eqref{eq:orbit}, and set $\gamma(P)=(\gamma_k(P))_{k=0}^{\lfloor(n-2)/2\rfloor}$.
By \eqref{eq:Ppartition} and the uniqueness of the expansion, these coefficients are independent of the choice of canonical labeling.

A \emph{maximal increasing run} of $\pi$ is a maximal increasing consecutive subsequence of $\pi$.
For example, the permutation $\pi=142365$ has three maximal increasing runs: $14$, $236$, and $5$.

\begin{lemma}\label{lem:runs}
Let $(P,\omega)$ be a canonically labeled rank-one poset, and let $\pi\in\mathcal L(P,\omega)$.
For every integer $k\ge0$, we have $\pi\in\mathcal R_k(P,\omega)$ if and only if $\pi$ has exactly $k+1$ maximal increasing runs, each of length at least two.
\end{lemma}

\begin{proof}
The descents of $\pi$ occur precisely between consecutive maximal increasing runs. Hence $\des(\pi)=k$ if and only if $\pi$ has exactly $k+1$ maximal increasing runs.

For $1<i<n$, the entry $\pi_i$ forms a maximal increasing run of length one if and only if $\pi_{i-1}>\pi_i>\pi_{i+1}$, which is the double-descent condition at position $i$.

Since $\pi$ is the label word of a linear extension, $\pi_1$ labels a minimal element and $\pi_n$ labels a maximal element. 
Since the labeling is canonical, $\pi_1<0<\pi_n$. 
Thus the first run has length one if and only if $0>\pi_1>\pi_2$, and the last run has length one if and only if $\pi_{n-1}>\pi_n>0$.
Under the convention $\pi_0=\pi_{n+1}=0$, these are the double-descent conditions at positions $1$ and $n$, respectively.

Therefore $0\pi0$ has no double descents if and only if every maximal increasing run of $\pi$ has length at least two. 
The two equivalences give the result by the definition of $\mathcal R_k(P,\omega)$. 
This completes the proof.
\end{proof}

An \emph{order ideal} of $P$ is a subset $I\subseteq P$ such that $y\in I$ and $x\le_P y$ imply $x\in I$.
The following result follows directly from~\cite[Proposition 3.5.2]{StanleyEC1}.

\begin{lemma}\label{lem:initial-segment}
Let $x_1,\ldots,x_n$ be a linear extension of $P$.
Then $\{x_1,\ldots,x_s\}$ is an order ideal of $P$ for every $1\le s\le n$.
\end{lemma}

For the remainder of this section, fix a canonically labeled rank-one poset $(P,\omega)$ on $n$ elements.
For each nonempty subset $A\subseteq P$, write $\omega(A)=\{\omega(x):x\in A\}$.
The vertices of $\Delta(P,\omega)$ are the order ideals $I$ of $P$ satisfying
\begin{equation}\label{eq:vertex}
2\le |I|\le n-2,\qquad
\max\omega(I)>\min\omega(P\setminus I).
\end{equation}

For a chain of vertices $I_1\subsetneq\cdots\subsetneq I_k$, set $I_0=\varnothing$ and $I_{k+1}=P$. 
Define $B_j=I_{j+1}\setminus I_j$ for $0\le j\le k$.
The blocks $B_0,\ldots,B_k$ are pairwise disjoint and their union is $P$. 
The block conditions are
\begin{equation}\label{eq:blocks}
\begin{aligned}
|B_j| &\ge 2 &&\text{for } 0\le j\le k,\\
\max\omega(B_j) &> \min\omega(B_{j+1})
&&\text{for } 0\le j<k.
\end{aligned}
\end{equation}
A nonempty set of $k$ vertices is a face if the vertices form such a chain and the associated blocks satisfy~\eqref{eq:blocks}.
The empty set is also a face.
We define $\Delta(P,\omega)$ to be the family of all these faces. 

\begin{lemma}\label{lem:complex}
The family $\Delta(P,\omega)$ is a simplicial complex.
\end{lemma}

\begin{proof}
For every vertex $I$, condition~\eqref{eq:vertex} ensures that the two sets $I$ and $P\setminus I$ satisfy~\eqref{eq:blocks} for $k=1$.
Hence $\{I\}$ is a face.

It remains to show that every subset of a face is again a face.
It suffices to verify that deleting any single vertex from a nonempty face produces another face. 

Let $F=\{I_1,\ldots,I_k\}$ be a nonempty face, with the vertices indexed so that
\[
I_1\subsetneq\cdots\subsetneq I_k.
\]
If $k=1$, deleting the only vertex produces the empty set, which is a face by definition.
We may therefore assume $k\ge2$.

Set $I_0=\varnothing$ and $I_{k+1}=P$.
The blocks determined by this chain are $B_i=I_{i+1}\setminus I_i$ for $0\le i\le k$.
Fix $j$ with $1\le j\le k$ and delete $I_j$.
The remaining vertices still form a strictly increasing chain of order ideals.
The blocks determined by this chain are unchanged except that $B_{j-1}$ and $B_j$ are replaced by the single set
\[
C=I_{j+1}\setminus I_{j-1}
=B_{j-1}\cup B_j.
\]

We check the two conditions in~\eqref{eq:blocks} for these blocks.
Since $B_{j-1}$ and $B_j$ are disjoint and each contains at least two elements, we have $|C|=|B_{j-1}|+|B_j|\ge4.$
All other blocks remain unchanged and contain at least two elements.
Thus the first condition in~\eqref{eq:blocks} holds.

For the second condition in~\eqref{eq:blocks}, only pairs of consecutive blocks involving $C$ need to be checked.
If $j=1$, then $C=B_0\cup B_1$ is the first block, and the next block is $B_2$.
Since $B_1\subseteq C$, condition~\eqref{eq:blocks} implies
$\max\omega(C) \ge \max\omega(B_1)> \min\omega(B_2).$

If $2\le j\le k-1$, then the blocks directly before and after $C$ are $B_{j-2}$ and $B_{j+1}$, respectively.
Using~\eqref{eq:blocks} and the inclusions $B_{j-1}\subseteq C$ and $B_j\subseteq C$, we obtain
$\max\omega(B_{j-2})>
\min\omega(B_{j-1})
\ge
\min\omega(C)$ and 
$\max\omega(C)
\ge
\max\omega(B_j)>
\min\omega(B_{j+1}).$

If $j=k$, then $C=B_{k-1}\cup B_k$ is the last block, and the preceding block is $B_{k-2}$.
Since $B_{k-1}\subseteq C$, condition~\eqref{eq:blocks} implies
$\max\omega(B_{k-2})>
\min\omega(B_{k-1})
\ge
\min\omega(C).$

Every other pair of consecutive blocks is unchanged.
Thus the second condition in~\eqref{eq:blocks} holds after deleting $I_j$.
It follows that $F\setminus\{I_j\}$ is a face.
By repeatedly deleting vertices, every subset of $F$ is again a face.
Thus $\Delta(P,\omega)$ is a simplicial complex.
\end{proof}

Our construction is related to the complex on linear extensions described by Stanley; see \cite[Theorem~3.13.1 and Exercise~3.152(a)]{StanleyEC1}.
The following bijection uses the classical correspondence between chains of order ideals and linear extensions; see~\cite[Theorem~3.13.1]{StanleyEC1}.
\begin{proposition}\label{prop:bijection}
For each integer $k\ge0$, there is a bijection $\phi$ between  the $k$-vertex faces of $\Delta(P,\omega)$ and the set $\mathcal R_k(P,\omega)$.
\end{proposition}

\begin{proof}
For a $k$-vertex face $F=\{I_1,\ldots,I_k\}$, index the order ideals by inclusion.
Set $I_0=\varnothing$, $I_{k+1}=P$, and $B_i=I_{i+1}\setminus I_i$ for $0\le i\le k$.
When $k=0$, the only block is $B_0=P$, of size $n\ge2$. 
For $0\le i\le k$, write
\begin{align*}
B_i=\{x_{i,1},\ldots,x_{i,r_i}\},
\qquad r_i=|B_i|,
\end{align*}
where the elements are indexed so that $\omega(x_{i,1})<\cdots<\omega(x_{i,r_i})$.
Let $w_i=\omega(x_{i,1})\omega(x_{i,2})\cdots\omega(x_{i,r_i})$.
Concatenating the words in the order $w_0,w_1,\ldots,w_k$, define
\begin{align*}
\phi(F)=w_0w_1\cdots w_k.
\end{align*}
Since the blocks partition $P$ and $\omega$ is injective, every label $\omega(x)$ with $x\in P$ occurs exactly once in $\phi(F)$.

We verify that $\phi(F)\in\mathcal L(P,\omega)$.
Suppose $x<_P y$ but $\omega(y)$ precedes $\omega(x)$ in $\phi(F)$.
Since $\omega$ is natural and every $w_i$ is increasing, $x$ and $y$ cannot belong to the same block.
Hence $y\in B_s$ and $x\in B_t$ for some $s<t$.
Then $I_{s+1}=B_0\cup\cdots\cup B_s$ contains $y$ but not $x$, contradicting the order-ideal property of $I_{s+1}$.
Therefore $\phi(F)\in\mathcal L(P,\omega)$.

By~\eqref{eq:blocks}, we have $\max\omega(B_i)>\min\omega(B_{i+1})$  for $0\le i<k$.
Since every $w_i$ is increasing, the words $w_0,\ldots,w_k$ are exactly the maximal increasing runs of $\phi(F)$.
Each run has length $|B_i|\ge2$.
By Lemma~\ref{lem:runs}, $\phi(F)\in\mathcal R_k(P,\omega)$.

To construct $\phi^{-1}$, let $\pi\in\mathcal R_k(P,\omega)$.
By Lemma~\ref{lem:runs}, write $\pi=w_0\cdots w_k$ as the decomposition into maximal increasing runs, where $w_i=\omega(x_{i,1})\cdots\omega(x_{i,r_i})$, $\omega(x_{i,1})<\cdots<\omega(x_{i,r_i})$ and $ r_i\ge2$. 
Let $B_i=\{x_{i,1},\ldots,x_{i,r_i}\}$ for $0\le i\le k$, and let $I_j=B_0\cup\cdots\cup B_{j-1}$ for $1\le j\le k$.
Define
\begin{align*}
\phi^{-1}(\pi)=\{I_1,\ldots,I_k\},
\end{align*}
with $\phi^{-1}(\pi)=\varnothing$ when $k=0$.

By Lemma~\ref{lem:initial-segment}, every $I_j$ is an order ideal.
Since the blocks are nonempty and pairwise disjoint, $I_1\subsetneq\cdots\subsetneq I_k$.
Moreover, since $|B_i|\ge2$, we have $|I_j|\ge2j\ge2$ and  $|P\setminus I_j|\ge2(k+1-j)\ge2$ for $1\le j\le k$.

Since consecutive runs are separated by descents, we have $\max\omega(B_{j-1})>\min\omega(B_j)$ for $1\le j\le k$.
As $B_{j-1}\subseteq I_j$ and $B_j\subseteq P\setminus I_j$, we obtain
\begin{align*}
\max\omega(I_j)
\ge \max\omega(B_{j-1})
> \min\omega(B_j)
\ge \min\omega(P\setminus I_j).
\end{align*}
Since $2\le |I_j|\le n-2$, every $I_j$ satisfies~\eqref{eq:vertex}.

With $I_0=\varnothing$ and $I_{k+1}=P$, the successive differences $I_{i+1}\setminus I_i$ are precisely $B_i$.
By the run lengths and the descents between consecutive runs, condition~\eqref{eq:blocks} holds.
Thus $\phi^{-1}(\pi)$ is a $k$-vertex face.

Cutting $\phi(F)$ into maximal increasing runs recovers the labels of the original blocks.
Since $\omega$ is injective, taking the corresponding element sets and their prefix unions recovers $F$.
Therefore $\phi^{-1}(\phi(F))=F$.
Conversely, applying $\phi$ to $\phi^{-1}(\pi)$ lists the labels in each block $B_i$ in increasing order, recovering the run $w_i$. 
The blocks are read in the order $B_0,\ldots,B_k$, so $\phi(\phi^{-1}(\pi))=w_0\cdots w_k=\pi$.
Hence $\phi$ is the desired bijection.
\end{proof}

Propositions~\ref{prop:gamma-count} and~\ref{prop:bijection} imply the following face-number identity.
\begin{equation}\label{eq:faces}
\sum_{\pi\in\Lext(P,\omega)}t^{\des(\pi)}
=\sum_{k=0}^{\lfloor(n-2)/2\rfloor}f_{k-1}\bigl(\Delta(P,\omega)\bigr)
t^k(1+t)^{n-2-2k}.
\end{equation}

\begin{corollary}\label{cor:face-gamma}
Let $(P,\omega)$ be a canonically labeled rank-one poset on $n$ elements. 
For every integer $k$ with $0\leq k\leq \lfloor(n-2)/2\rfloor$, we have
\begin{align*}
f_{k-1}\bigl(\Delta(P,\omega)\bigr)
=|\R_k(P,\omega)|
=\gamma_k(P).
\end{align*}
\end{corollary}

\section{Dimension and balancedness}
In this section, we determine the dimension of $\Delta(P,\omega)$ and study its balancedness.
\begin{lemma}\label{lem:colors}
For every canonically labeled rank-one poset $(P,\omega)$ on $n\geq 2$ elements, the map
\[\operatorname{col}(I)=\left\lfloor\frac{|I|}{2}\right\rfloor\] 
is a proper coloring of $\Delta(P,\omega)$ with colors in $\{1,\ldots,\lfloor(n-2)/2\rfloor\}$. 
Moreover, every face has at most $\lfloor(n-2)/2\rfloor$ vertices, and hence $\dim\Delta(P,\omega)\leq \lfloor(n-2)/2\rfloor-1$.
\end{lemma}

\begin{proof}
By~\eqref{eq:vertex}, every vertex $I$ satisfies $2\leq |I|\leq n-2$. 
Hence $1\leq\operatorname{col}(I)\leq\lfloor(n-2)/2\rfloor$.
Let $I_j$ and $I_{j+1}$ be two consecutive vertices in a face. 
By \eqref{eq:blocks}, we obtain $|I_{j+1}|-|I_j|=|I_{j+1}\setminus I_j|=|B_j|\geq 2$.
It follows that $\operatorname{col}(I_j)<\operatorname{col}(I_{j+1})$.
For every nonempty face $\{I_1,\ldots,I_k\}$ with $I_1\subsetneq\cdots\subsetneq I_k$, we have $\operatorname{col}(I_1)<\cdots<\operatorname{col}(I_k)$.
Hence the vertices of every face have pairwise distinct colors, so the coloring is proper.

Finally, since the coloring uses at most $\lfloor(n-2)/2\rfloor$ colors, every face has at most $\lfloor(n-2)/2\rfloor$ vertices. 
Consequently, $\dim\Delta(P,\omega)\leq \lfloor(n-2)/2\rfloor-1$.
\end{proof}

Let $(P,\omega)$ be a canonically labeled rank-one poset on $n$ elements. 
Given a linear extension $(x_1,\ldots,x_n)$ of $P$, choose $k\ge0$ descent positions $s_1<\cdots<s_k$ of the label word $\omega(x_1)\cdots\omega(x_n)$.
Set $s_0=0$, $s_{k+1}=n$, and $I_j=\{x_1,\ldots,x_{s_j}\}$ for $1\le j\le k$.

\begin{lemma}\label{lem:spaced-descents}
If $s_{j+1}-s_j\ge2$ for $0\le j\le k$, then $\{I_1,\ldots,I_k\}$ is a $k$-vertex face of $\Delta(P,\omega)$.
\end{lemma}

\begin{proof}
For $k=0$, the collection is the empty face.
Assume $k\ge1$. By Lemma~\ref{lem:initial-segment}, each $I_j$ is an order ideal. The strict inequalities among the $s_j$ give $I_1\subsetneq\cdots\subsetneq I_k$.

With $I_0=\varnothing$ and $I_{k+1}=P$, the associated blocks are
\begin{align*}
B_j=I_{j+1}\setminus I_j
=\{x_{s_j+1},\ldots,x_{s_{j+1}}\}
\qquad\text{for }0\le j\le k.
\end{align*}
For $0\le j\le k$, we have $|B_j|=s_{j+1}-s_j\ge2$.
Since $s_0=0$ and $s_{k+1}=n$, we have $s_1\ge2$ and $s_k\le n-2$. Hence $2\le |I_j|=s_j\le n-2$ for $1\le j\le k$.

Since each $s_j$ is a descent position, for $1\le j\le k$ we have
\[
\max\omega(I_j)
\ge \max\omega(B_{j-1})
\ge \omega(x_{s_j})
> \omega(x_{s_j+1})
\ge \min\omega(B_j)
\ge \min\omega(P\setminus I_j).
\]
Thus the ideals satisfy~\eqref{eq:vertex} and the blocks satisfy~\eqref{eq:blocks}, so $\{I_1,\ldots,I_k\}$ is a face.
\end{proof}

Define the sets of minimal and maximal elements of $P$ by
\begin{align*}
\Min(P) &= \{x\in P : \nexists\,y\in P \text{ with } y<_P x\},\\
\Max(P) &= \{x\in P : \nexists\,y\in P \text{ with } x<_P y\}.
\end{align*}
For a nonempty subset $A=\{z_1,\ldots,z_\ell\}\subseteq P$, choose the permutation $(i_1,\ldots,i_\ell)$ of $(1,\ldots,\ell)$ such that $\omega(z_{i_1})>\cdots>\omega(z_{i_\ell})$.
Define 
\begin{align*}
D(A)=(z_{i_1},\ldots,z_{i_\ell}).
\end{align*}  
For sequences $X=(x_1,\ldots,x_p)$ and $Y=(y_1,\ldots,y_q)$, define their concatenation by
\begin{align*}
XY=(x_1,\ldots,x_p,y_1,\ldots,y_q).
\end{align*}

Write $(x_1,\ldots,x_n)=D(\Min(P))\,D(\Max(P))$.
The sequence lists all minimal elements before all maximal elements. 
Since every strict relation in $P$ has the form $x<_P y$ with $x\in\Min(P)$ and $y\in\Max(P)$, the sequence is a linear extension of $P$.

For $1\leq i<n$ with $i\ne|\Min(P)|$, the elements $x_i$ and $x_{i+1}$ belong to the same rank level.
By the definition of $D$, we have $\omega(x_i)>\omega(x_{i+1})$.
At $i=|\Min(P)|$, the element $x_i$ is the last minimal element and $x_{i+1}$ is the first maximal element, so $\omega(x_i)<0<\omega(x_{i+1})$.
Thus the label word has a descent at every position except $|\Min(P)|$.

Set $E=\{2j:1\le j\le\lfloor(n-2)/2\rfloor\}$.
The indices in each of $E$ and $\{n-i:i\in E\}$ lie between $2$ and $n-2$, and any two distinct indices in the same set differ by at least two.
If $S$ is a subset of either set and every index in $S$ is a descent position of the label word of a linear extension of $P$, then Lemma~\ref{lem:spaced-descents} gives a face of $\Delta(P,\omega)$ with $|S|$ vertices.

\begin{lemma}\label{lem:dimension-exception}
Suppose $(P,\omega)$ is a canonically labeled rank-one poset on $n$ elements. 
If both rank levels have an even number of elements and every minimal element is below every maximal element, then $\dim\Delta(P,\omega)=\lfloor(n-2)/2\rfloor-2$.
\end{lemma}

\begin{proof}
Since the two rank levels are nonempty and each has an even number of elements, $n$ is even and $|\Min(P)|$ belongs to $E=\{2,4,\ldots,n-2\}$. 
In the linear extension $D(\Min(P))\,D(\Max(P))$, the only ascent is at position $|\Min(P)|$. 
Thus every index in $E\setminus\{|\Min(P)|\}$ is a descent. 
Applying Lemma~\ref{lem:spaced-descents} to these indices yields a face with $|E|-1=n/2-2$ vertices.

By Lemma~\ref{lem:colors}, a face has at most $n/2-1$ vertices. 
Suppose a face has $n/2-1$ vertices, with order ideals $I_1\subsetneq\cdots\subsetneq I_{n/2-1}$.
Set $I_0=\varnothing$, $I_{n/2}=P$, and $B_j=I_{j+1}\setminus I_j$ for $0\leq j<n/2$.
By~\eqref{eq:blocks}, each block has at least two elements.
The $n/2$ blocks are pairwise disjoint and together contain all $n$ elements of $P$, so
\begin{align*}
\sum_{j=0}^{n/2-1}|B_j|=n=2\cdot\frac{n}{2}.
\end{align*}
Since each of the $n/2$ summands is at least $2$ and their sum is $n$, we have $|B_j|=2$ for $0\leq j<n/2$.
Consequently,
\begin{align*}
|I_j|=\sum_{\ell=0}^{j-1}|B_\ell|=2j
\qquad\text{for }1\leq j<n/2.
\end{align*}
Since $|\Min(P)|\in E$, the index $j=|\Min(P)|/2$ lies in this range. 
Thus the face contains the ideal $I_{|\Min(P)|/2}$ of size $|\Min(P)|$.

Every maximal element lies above every minimal element.
Therefore, an order ideal containing a maximal element has at least $|\Min(P)|+1$ elements. It follows that the only order ideal of size $|\Min(P)|$ is $\Min(P)$.
But $P\setminus\Min(P)=\Max(P)$, and the canonical labeling satisfies
\begin{align*}
\max\omega(\Min(P))<0<\min\omega(\Max(P)).
\end{align*}
Thus $\Min(P)$ fails the vertex condition~\eqref{eq:vertex}, a contradiction.

The maximum number of vertices in a face is therefore $n/2-2$. 
Hence
\begin{align*}
\dim\Delta(P,\omega)
=\frac n2-3
=\left\lfloor\frac{n-2}{2}\right\rfloor-2.
\end{align*}
This completes the proof.
\end{proof}

\begin{proposition}\label{prop:dimension}
Suppose $(P,\omega)$ is a canonically labeled rank-one poset on $n\geq 2$ elements. 
If at least one rank level has an odd number of elements or some minimal element is incomparable with some maximal element, then $\dim\Delta(P,\omega)=\lfloor\frac{n-2}{2}\rfloor-1$.
\end{proposition}

\begin{proof}
By Lemma~\ref{lem:colors}, it suffices to construct a face with $\lfloor(n-2)/2\rfloor$ vertices.
We consider two cases according to whether both rank levels have an even number of elements.

\smallskip
\noindent
\textbf{Case 1. At least one rank level has an odd number of elements.}
Use the linear extension $D(\Min(P))\,D(\Max(P))$.
If $|\Min(P)|$ is odd, select the indices in $E$.
Every index in $E$ is even, so the selected indices avoid the ascent at $|\Min(P)|$. 
Otherwise, $|\Max(P)|$ is odd, so $n-|\Min(P)|=|\Max(P)|\notin E$. 
Select the indices in $\{n-i:i\in E\}$, which therefore excludes $|\Min(P)|$.
All selected indices are descent positions.
By Lemma~\ref{lem:spaced-descents}, they determine a face with $|E|=\lfloor(n-2)/2\rfloor$ vertices.

\smallskip
\noindent
\textbf{Case 2. Both rank levels have an even number of elements.}
In this case, the proposition assumes that some minimal element is incomparable with some maximal element. 
Choose such a pair $v\in\Min(P)$ and $u\in\Max(P)$, and form the concatenation
\begin{align*}
(x_1,\ldots,x_n)
=D(\Min(P)\setminus\{v\})\,(u,v)\,
D(\Max(P)\setminus\{u\}).
\end{align*}
If $x<_P u$, then $x\in\Min(P)\setminus\{v\}$, so $x$ appears before $u$. 
Every maximal element other than $u$ appears after all minimal elements. 
Thus $x_1,\ldots,x_n$ is a linear extension of $P$.

The first sequence has $|\Min(P)|-1$ elements, so $x_{|\Min(P)|}=u$ and $x_{|\Min(P)|+1}=v$. 
Hence
\begin{align*}
\omega(x_{|\Min(P)|})=\omega(u)>0>\omega(v)
=\omega(x_{|\Min(P)|+1}).
\end{align*}
The labels decrease within each of $D(\Min(P)\setminus\{v\})$ and $D(\Max(P)\setminus\{u\})$.
Since $|\Min(P)|$ is even, the two remaining indices $|\Min(P)|-1$ and $|\Min(P)|+1$ are odd. Every index in $E$ is even, so
\begin{align*}
E\cap\{|\Min(P)|-1,\ |\Min(P)|+1\}=\varnothing.
\end{align*}
Therefore every index in $E$ is a descent position.
By Lemma~\ref{lem:spaced-descents}, the indices in $E$ determine a face with $|E|=\lfloor(n-2)/2\rfloor$ vertices.
This completes the proof.
\end{proof}

\begin{proof}[\textbf{Proof of Theorem~\ref{thm:general}}]
By Lemma~\ref{lem:complex}, $\Delta(P,\omega)$ is a simplicial complex. 
The face-number identity follows from~\eqref{eq:faces}, and Lemma~\ref{lem:colors} provides the proper coloring.
The dimension formula follows from Lemma~\ref{lem:dimension-exception} and Proposition~\ref{prop:dimension}, whose hypotheses together cover all canonically labeled rank-one posets.
This completes the proof.
\end{proof}

\begin{lemma}[{\cite[Theorem~1.2(i), p.~170]{FFK}}]
\label{lem:ffk}
If a finite simplicial complex is properly colorable with $r\geq 1$ colors, then its face vector satisfies the Frankl--F\"uredi--Kalai inequalities for $r$-colorable complexes.
\end{lemma}

By Theorem~\ref{thm:general} and Lemma~\ref{lem:ffk}, we obtain the following result.

\begin{corollary}\label{cor:ffk}
For a rank-one poset $P$ on $n\geq 4$ elements, the $\gamma$-vector $\gamma(P)$ satisfies the Frankl--F\"uredi--Kalai inequalities for $\lfloor(n-2)/2\rfloor$-colorable simplicial complexes.
\end{corollary}

\begin{proof}
Choose a canonical labeling $\omega$ of $P$.
By Lemma~\ref{lem:colors}, $\Delta(P,\omega)$ is $\lfloor(n-2)/2\rfloor$-colorable. 
Thus, by Lemma~\ref{lem:ffk}, the face numbers of $\Delta(P,\omega)$ satisfy the Frankl--F\"uredi--Kalai inequalities for $\lfloor(n-2)/2\rfloor$-colorable complexes.
By Corollary~\ref{cor:face-gamma}, $\gamma_k(P)=f_{k-1}(\Delta(P,\omega))$ for $0\le k\le\lfloor(n-2)/2\rfloor$, so the same inequalities hold for $\gamma(P)$. 
This completes the proof.
\end{proof}

\begin{proof}[\textbf{Proof of Theorem~\ref{cor:fence}}]
For $n=1$, we have $U_1(t)=1$ and $\gamma_1=(1)$.
Take $\Delta_1=\{\varnothing\}$. 
This complex is balanced of dimension $-1$, and its $f$-vector is $(1)$.
The Frankl--F\"uredi--Kalai inequalities are vacuous in this case.

Assume now that $n\ge2$.
Recall that $P_n$ is the zig-zag poset on $p_1,\ldots,p_n$. 
Its cover relations are $p_{2j-1}<_P p_{2j}$ for $1\leq j\leq\lfloor n/2\rfloor$ and $p_{2j+1}<_P p_{2j}$ for $1\leq j\leq\lfloor(n-1)/2\rfloor$.

Define the labeling $\omega$ by
\begin{align*}
\omega(p_{2i-1})=i-\lceil n/2\rceil-1,\qquad
\omega(p_{2j})=j
\end{align*}
for $1\leq i\leq\lceil n/2\rceil$ and $1\leq j\leq\lfloor n/2\rfloor$.
The odd-indexed elements are minimal and receive distinct negative labels. 
The even-indexed elements are maximal and receive distinct positive labels.
Thus $\omega$ is a canonical labeling.

Set $\Delta_n=\Delta(P_n,\omega)$.
Since $U_n(t)$ is the descent polynomial of $(P_n,\omega)$, we have $\gamma_{n,k}=\gamma_k(P_n)$.
By Corollary~\ref{cor:face-gamma},
\begin{align}\label{eq:fence-face-numbers}
\gamma_{n,k}=f_{k-1}(\Delta_n)
\end{align}
for $0\leq k\leq\lfloor(n-2)/2\rfloor$.

If $n=2$ or $3$, no order ideal satisfies $2\leq |I|\leq n-2$. 
Hence $\Delta_n=\{\varnothing\}$ is balanced of dimension $-1$.
Since $\lfloor(n-2)/2\rfloor=0$, the only $\gamma$-coefficient is $\gamma_{n,0}=f_{-1}(\Delta_n)=1$.
The Frankl--F\"uredi--Kalai assertion is vacuous in these two cases.

Assume now that $n\geq4$. By the cover relations of $P_n$, the only element covering $p_1$ is $p_2$.
Since $P_n$ has rank one, $p_2$ is therefore the only element strictly above $p_1$.
Thus the minimal element $p_1$ and the maximal element $p_4$ are incomparable. 
By Proposition~\ref{prop:dimension},
\begin{align*}
\dim\Delta_n=\lfloor(n-2)/2\rfloor-1.
\end{align*}

By Lemma~\ref{lem:colors}, $\Delta_n$ admits a proper coloring with $\lfloor(n-2)/2\rfloor$ colors.
The number of colors is $\dim(\Delta_n)+1$, so $\Delta_n$ is balanced. 
By~\eqref{eq:fence-face-numbers}, the $\gamma$-vector of $U_n(t)$ is the face vector of the balanced simplicial complex $\Delta_n$.

Finally, Corollary~\ref{cor:ffk} applied to $P_n$ shows that the $\gamma$-vector of $U_n(t)$ satisfies the Frankl--F\"uredi--Kalai inequalities for $\lfloor(n-2)/2\rfloor$-colorable simplicial complexes.
This completes the proof.
\end{proof}

\section{Further applications}\label{sec:applications}

In this section, we give some further applications.

\begin{lemma}\label{prop:balanced-rank-one}
For every canonically labeled rank-one poset $(P,\omega)$, the simplicial complex $\Delta(P,\omega)$ constructed in Theorem~\ref{thm:general} is balanced.
\end{lemma}

\begin{proof}
By Theorem~\ref{thm:general}, it suffices to consider the case where both rank levels have even size and every minimal element is below every maximal element.
Let $n=|P|$ and $a=|\operatorname{Min}(P)|$.
By Lemma~\ref{lem:dimension-exception},
\begin{align*}
\dim\Delta(P,\omega)=\frac{n}{2}-3.
\end{align*}

We claim that the only order ideal of size $a$ is $\operatorname{Min}(P)$.
Indeed, an order ideal containing a maximal element must contain every minimal element, so it has at least $a+1$ elements. 
An order ideal of size $a$ therefore contains only minimal elements and must equal $\operatorname{Min}(P)$.
This proves the claim.

Since $\max\omega(\operatorname{Min}(P))<0<
\min\omega(\operatorname{Max}(P))$, this ideal fails condition~\eqref{eq:vertex}.
Thus every vertex $I$ of $\Delta(P,\omega)$ satisfies $|I|\neq a$.

Since no vertex has size $a$, define $\operatorname{col}'$ on the vertices of $\Delta(P,\omega)$ by
\begin{align*}
\operatorname{col}'(I)=
\begin{cases}
\displaystyle\left\lfloor\frac{|I|}{2}\right\rfloor,
& |I|<a,\\[6pt]
\displaystyle\left\lfloor\frac{|I|-1}{2}\right\rfloor,
& |I|>a.
\end{cases}
\end{align*}
By~\eqref{eq:vertex}, we have $1\le\operatorname{col}'(I)\le a/2-1$ when $|I|<a$, and $a/2\le\operatorname{col}'(I)\le n/2-2$ when $|I|>a$.

Let $I\subsetneq J$ be consecutive vertices in a face.
By~\eqref{eq:blocks}, we have $|J|-|I|\ge2$.
If $|I|,|J|<a$, then
\begin{align*}
\operatorname{col}'(J)
=\left\lfloor\frac{|J|}{2}\right\rfloor
\ge\left\lfloor\frac{|I|+2}{2}\right\rfloor
=\operatorname{col}'(I)+1.
\end{align*}
If $a<|I|,|J|$, then
\begin{align*}
\operatorname{col}'(J)
=\left\lfloor\frac{|J|-1}{2}\right\rfloor
\ge\left\lfloor\frac{|I|+1}{2}\right\rfloor
=\operatorname{col}'(I)+1.
\end{align*}
Finally, if $|I|<a<|J|$, the two color ranges give
\begin{align*}
\operatorname{col}'(I)
\le\frac a2-1
<\frac a2
\le\operatorname{col}'(J).
\end{align*}
Hence the colors strictly increase along the chain of ideals in each face. 
Thus $\operatorname{col}'$ is a proper coloring with at most $n/2-2$ colors.
As
\begin{align*}
\frac{n}{2}-2=\dim\Delta(P,\omega)+1,
\end{align*}
the complex $\Delta(P,\omega)$ is balanced.
This completes the proof.
\end{proof}

For a simplicial complex $\Delta$ of dimension $d-1$, define its $h$-polynomial by
\begin{align*}
h(\Delta;t)
=\sum_{i=0}^{d}
f_{i-1}(\Delta)t^i(1-t)^{d-i}.
\end{align*}
Reiner and Welker~\cite[Section~3.1]{ReinerWelker} constructed an equatorial sphere $\Delta_{\mathrm{eq}}(P)$ for a naturally labeled graded poset $(P,\omega)$.
They proved~\cite[Corollary~3.8]{ReinerWelker} that
\begin{align}\label{eq:equatorial-h}
h(\Delta_{\mathrm{eq}}(P);t)
=
\sum_{\pi\in\mathcal L(P,\omega)}
t^{\operatorname{des}(\pi)}.
\end{align}

Br\"and\'en~\cite[Section~9]{Branden} used this identity to give a combinatorial interpretation of the $\gamma$-coefficients of equatorial spheres.
For rank-one posets, we obtain the following result.
\begin{corollary}\label{prop:equatorial}
For a canonically labeled rank-one poset $(P,\omega)$ on $n\ge2$ elements,
\begin{align*}
\gamma_k(\Delta_{\mathrm{eq}}(P))
=
f_{k-1}(\Delta(P,\omega))
\end{align*}
for $0\le k\le\lfloor(n-2)/2\rfloor$.
Here $f_{k-1}(\Delta(P,\omega))=0$ when $k-1>\dim\Delta(P,\omega)$.
The complex $\Delta(P,\omega)$ is balanced, and the $\gamma$-vector of $\Delta_{\mathrm{eq}}(P)$ satisfies the Frankl--F\"uredi--Kalai inequalities.
\end{corollary}

\begin{proof}
Every maximal chain of $P$ has two elements, so $P$ is graded with two rank levels. 
The canonical labeling $\omega$ is natural. Replacing the labels by $1,\ldots,n$ in increasing order preserves the natural labeling condition and all descent positions. 
Thus \eqref{eq:equatorial-h} holds for $(P,\omega)$.

Set $d=n-2$. Since $P$ has two rank levels, the dimension formula in~\cite[Proposition~3.11]{ReinerWelker} gives $\dim\Delta_{\mathrm{eq}}(P)=n-3=d-1$.
For $n=2$, the equatorial complex is $\{\varnothing\}$, viewed as the $(-1)$-sphere.

By~\eqref{eq:equatorial-h} and~\eqref{eq:faces}, we obtain the $\gamma$-expansion of $h(\Delta_{\mathrm{eq}}(P);t)$ with parameter $d=\dim\Delta_{\mathrm{eq}}(P)+1$:
\begin{align*}
h(\Delta_{\mathrm{eq}}(P);t)
=\sum_{\pi\in\mathcal L(P,\omega)}
t^{\operatorname{des}(\pi)}
=\sum_{k=0}^{\lfloor d/2\rfloor}
f_{k-1}(\Delta(P,\omega))\,t^k(1+t)^{d-2k}.
\end{align*}
By the uniqueness of the $\gamma$-expansion and Corollary~\ref{cor:face-gamma}, we have
\begin{align*}
\gamma_k(\Delta_{\mathrm{eq}}(P))
=f_{k-1}(\Delta(P,\omega))
=\gamma_k(P)
\end{align*}
for $0\le k\le\lfloor d/2\rfloor$.

By Lemma~\ref{prop:balanced-rank-one}, $\Delta(P,\omega)$ is balanced. 
For $n=2$ or $3$, we have $\lfloor d/2\rfloor=0$ and $\gamma_0(\Delta_{\mathrm{eq}}(P))=f_{-1}(\Delta(P,\omega))=1$, so the Frankl--F\"uredi--Kalai inequalities hold.
For $n\ge4$, the $\gamma$-vector of $\Delta_{\mathrm{eq}}(P)$ equals $\gamma(P)$, which satisfies the Frankl--F\"uredi--Kalai inequalities by Corollary~\ref{cor:ffk}.
This completes the proof.
\end{proof}

Nevo and Petersen~\cite[Conjecture~6.3]{NP} conjectured that the $\gamma$-vector of every flag homology sphere satisfies the Frankl--F\"uredi--Kalai inequalities. 
Corollary~\ref{prop:equatorial} proves the conjecture for flag equatorial spheres of rank-one posets. 
The corollary also applies to the equatorial sphere of the six-element zig-zag poset, which Reiner and Welker~\cite[Example~3.22]{ReinerWelker} showed to be nonflag.
Athanasiadis~\cite[Theorem~3.1]{Athanasiadis2026} established the $\gamma$-positivity of the $h$-polynomials associated with integer compositions by interpreting them as Eulerian polynomials of graded posets of width two.

Stanley~\cite{Stanley1986} introduced the order polytope $\mathcal O(P)$ and the chain polytope $\mathcal C(P)$ of a finite poset $P$ and showed that they have the same Ehrhart polynomial. 
If $|P|=n$, both are $n$-dimensional lattice polytopes in $\mathbb R^P\cong\mathbb R^n$ \cite[Sections~1 and~2]{Stanley1986}.

For an $n$-dimensional lattice polytope $Q\subseteq\mathbb R^n$, the Ehrhart $h^*$-polynomial $h^*(Q;t)$ is defined by
\begin{align}\label{eq:ehrhart-hstar}
\sum_{q\ge0}|qQ\cap\mathbb Z^n|t^q
=\frac{h^*(Q;t)}{(1-t)^{n+1}}.
\end{align}

\begin{lemma}[{\cite[Theorem~4.1]{Stanley1986}}]
	\label{lem:poset-lattice-points}
	Let $P$ be a finite poset on $n$ elements.
	For every nonnegative integer $q$, we have
	\[
	|q\mathcal O(P)\cap\mathbb Z^n|
	=
	|q\mathcal C(P)\cap\mathbb Z^n|
	=
	\Omega_P(q+1).
	\]
\end{lemma}
For $q=0$, both dilates consist only of the origin, and $\Omega_P(1)=1$.
\begin{corollary}\label{prop:poset-polytopes}
Let $(P,\omega)$ be a canonically labeled rank-one poset on $n\ge2$ elements. 
Then
\begin{align}\label{eq:poset-polytopes-gamma}
\gamma_k(h^*(\mathcal O(P);t))
=\gamma_k(h^*(\mathcal C(P);t))
=f_{k-1}(\Delta(P,\omega))
\end{align}
for $0\le k\le\lfloor(n-2)/2\rfloor$, where face numbers above the dimension of $\Delta(P,\omega)$ are zero. 
In particular, these $\gamma$-vectors satisfy the Frankl--F\"uredi--Kalai inequalities.
\end{corollary}

\begin{proof}
By Lemma~\ref{lem:poset-lattice-points}
and~\eqref{eq:ehrhart-hstar}, we obtain
\begin{align*}
\frac{h^*(\mathcal O(P);t)}{(1-t)^{n+1}}
=\frac{h^*(\mathcal C(P);t)}{(1-t)^{n+1}}
=\sum_{q\ge0}\Omega_P(q+1)t^q
=\frac{\displaystyle
\sum_{\pi\in\mathcal L(P,\omega)}t^{\des(\pi)}}
{(1-t)^{n+1}}.
\end{align*}
The last equality follows by dividing both sides of~\eqref{eq:Ppartition} by $t$ and shifting the summation index to start at $q=0$.

Multiplying by $(1-t)^{n+1}$ and applying~\eqref{eq:faces}, we obtain
\begin{align*}
h^*(\mathcal O(P);t)
=h^*(\mathcal C(P);t)
=\sum_{\pi\in\mathcal L(P,\omega)}t^{\des(\pi)}
=\sum_{k=0}^{\lfloor(n-2)/2\rfloor}
f_{k-1}(\Delta(P,\omega))
t^k(1+t)^{n-2-2k}.
\end{align*}
The $k=0$ term is $(1+t)^{n-2}$, since $f_{-1}(\Delta(P,\omega))=1$. 
For $k\ge1$, the basis polynomial $t^k(1+t)^{n-2-2k}$ has degree $n-2-k<n-2$.
Thus the common $h^*$-polynomial has degree $n-2$.
The uniqueness of the $\gamma$-expansion yields~\eqref{eq:poset-polytopes-gamma}.

By Lemma~\ref{prop:balanced-rank-one}, $\Delta(P,\omega)$ is balanced.
For $n=2$ or $3$, the common $\gamma$-vector is $(1)$, so the Frankl--F\"uredi--Kalai inequalities hold.
For $n\ge4$, by~\eqref{eq:poset-polytopes-gamma} and Corollary~\ref{cor:face-gamma}, the common $\gamma$-vector equals $\gamma(P)$. 
The Frankl--F\"uredi--Kalai inequalities then follow from Corollary~\ref{cor:ffk}.
This completes the proof.
\end{proof}

For an integer $m\ge3$, let $C_{2m}$ be the crown poset obtained from the zig-zag poset $P_{2m}$ by adding the cover relation $p_1<p_{2m}$. 
Crown posets are special cases of circular fence posets, whose polyhedral properties were studied by Kantarc{\i} O\u{g}uz, \"Ozel, and Ravichandran~\cite{KOR2024}.
Lundstr\"om and Saud Maia Leite~\cite[Theorem~5.19]{LundstromLeite} proved that the Ehrhart $h^*$-polynomial of $\mathcal O(C_{2m})$ is $\gamma$-nonnegative.
The following corollary realizes its $\gamma$-coefficients as face numbers of a balanced simplicial complex.

\begin{corollary}
For every $m\ge3$ and every canonical labeling $\omega$ of $C_{2m}$, the complex $\Delta(C_{2m},\omega)$ is balanced of dimension $m-2$, and
\begin{align}\label{eq:crown-faces}
\gamma_k(h^*(\mathcal O(C_{2m});t))
=f_{k-1}(\Delta(C_{2m},\omega))
\end{align}
for $0\le k\le m-1$.
Consequently, the $\gamma$-vector of $h^*(\mathcal O(C_{2m});t)$ satisfies the Frankl--F\"uredi--Kalai inequalities.
\end{corollary}

\begin{proof}
Fix an arbitrary canonical labeling $\omega$ of $C_{2m}$.
The added cover runs from the minimal element $p_1$ to the maximal element $p_{2m}$ of $P_{2m}$.
Thus $C_{2m}$ is a rank-one poset with $m$ elements in each rank level. In $P_{2m}$, the minimal element $p_1$ has one upper cover, and every other minimal element has two. 
After adding $p_1<p_{2m}$, each minimal element lies below exactly two maximal elements.
Since $m\ge3$, there is an incomparable pair consisting of a minimal element and a maximal element.
As $|C_{2m}|=2m$, Theorem~\ref{thm:general} gives
\begin{align*}
\dim\Delta(C_{2m},\omega)
=\left\lfloor\frac{2m-2}{2}\right\rfloor-1
=m-2.
\end{align*}

By Lemma~\ref{prop:balanced-rank-one}, $\Delta(C_{2m},\omega)$ is balanced.
Applying Corollary~\ref{prop:poset-polytopes} to $P=C_{2m}$ gives~\eqref{eq:crown-faces} and the Frankl--F\"uredi--Kalai inequalities.
This completes the proof.
\end{proof}

\section*{Declaration of Generative AI}
The authors used generative AI to assist with research and writing.
They verified all arguments and results and take full
responsibility for the manuscript.


\begin{thebibliography}{99}
\bibitem{Alex}
P.~Alexandersson,
Real-rooted Eulerian polynomials from permutations, words, and paths,
preprint (2026), available at: 
arXiv:2609.07325v2.

\bibitem{AlexJal}
P.~Alexandersson and A.~Jal,
Rook matroids and log-concavity of $P$-Eulerian polynomials,
\emph{S\'em. Lothar. Combin.}
\textbf{93B} (2025), Article~59, 12 pp.

\bibitem{Athanasiadis2026}
C.~A.~Athanasiadis,
Lattice point enumeration of polytopes associated to integer compositions,
\emph{Ann. Comb.} (2026).

\bibitem{Branden}
P.~Br\"and\'en,
Actions on permutations and unimodality of descent polynomials,
\emph{European J. Combin.} \textbf{29} (2008), no.~2, 514--531.

\bibitem{FFK}
P.~Frankl, Z.~F\"uredi, and G.~Kalai,
Shadows of colored complexes,
\emph{Math. Scand.} \textbf{63} (1988), 169--178.

\bibitem{KOR2024}
E.~Kantarc{\i} O\u{g}uz, C.~Y.~\"Ozel
and M.~Ravichandran,
Chainlink polytopes and Ehrhart equivalence,
\emph{Ann. Comb.} \textbf{28} (2024), 1141--1166.

\bibitem{LundstromLeite}
T.~Lundstr\"om and L.~Saud Maia Leite,
Order polytopes of crown posets,
\emph{European J. Combin.} \textbf{133} (2026), 104304.


\bibitem{NPT}
E.~Nevo, T.~K.~Petersen, and B.~E.~Tenner,
The $\gamma$-vector of a barycentric subdivision,
\emph{J. Combin. Theory Ser. A} \textbf{118} (2011), no.~4, 1364--1380.

\bibitem{NP}
E.~Nevo and T.~K.~Petersen,
On $\gamma$-vectors satisfying the Kruskal--Katona inequalities,
\emph{Discrete Comput. Geom.} \textbf{45} (2011), 503--521.

\bibitem{PZ}
T.~K.~Petersen and Y.~Zhuang,
Zig-zag Eulerian polynomials,
\emph{European J. Combin.} \textbf{124} (2025),~104073.

\bibitem{ReinerWelker}
V.~Reiner and V.~Welker,
On the Charney--Davis and Neggers--Stanley conjectures,
\emph{J. Combin. Theory Ser. A}
\textbf{109} (2005), 247--280.

\bibitem{Stanley1986}
R.~P.~Stanley,
Two poset polytopes,
\emph{Discrete Comput. Geom.}
\textbf{1} (1986), 9--23.

\bibitem{StanleyEC1}
R.~P.~Stanley,
\emph{Enumerative Combinatorics, Volume~1},
second edition,
Cambridge University Press, Cambridge, 2012.

\bibitem{Torres}
E.~Torres,
Real-rootedness and log-concavity of zigzag Eulerian polynomials,
preprint (2026), available at: arXiv:2611.22110v1.


\end{thebibliography}
\end{document}